\documentclass[12pt]{article} 
\usepackage[reqno]{amsmath} % equation numbers on the right
\usepackage{amssymb, amsthm}
\usepackage{enumerate}  % lets you do \begin{enumerate}[(a)] and such
\usepackage{url} % formats url's nicely using \url{}
\usepackage{hyperref}

\newtheoremstyle{plainsl}%
	{\topsep}
	{\topsep}
	{\slshape} % only non-default setting
	{}
	{\normalfont\bfseries}
	{.}
	{ }
	{}

\swapnumbers

\theoremstyle{plainsl}
\newtheorem{theorem}{Theorem}[section]
\newtheorem{lemma}[theorem]{Lemma}
\newtheorem{corollary}[theorem]{Corollary}

\newcommand\cref[1]{Corollary~\ref{cor:#1}}

\renewcommand\proof{\noindent\textsl{Proof. }}
\newcommand\sqr[2]{{\vbox{\hrule height.#2pt
    \hbox{\vrule width.#2pt height#1pt \kern#1pt
        \vrule width.#2pt}\hrule height.#2pt}}}
\renewcommand\qed{%
	\ifmmode\eqno\sqr53
	\else\nolinebreak\ \hfill\sqr53\medbreak\fi}

\numberwithin{equation}{section}

\newcommand{\ii}{\mathrm{i}}
\newcommand{\e}{\mathrm{e}}
\newcommand{\ee}{\mathbf{e}}

\title{Perfect state transfer in products and covers of graphs}
\author{G.~Coutinho\thanks{Department of Combinatorics and Optimization, University of Waterloo, \texttt{ \{gcoutinho, cgodsil\}@uwaterloo.ca}}\ \thanks{Supported by Capes Foundation, Ministry of Education, Brazil.} , C.~Godsil\footnotemark[1]}

\begin{document}

\maketitle

\begin{abstract}
A continuous-time quantum walk on a graph $X$ is represented by the complex matrix $\exp (-\ii t A)$, where $A$ is the adjacency matrix of $X$ and $t$ is a non-negative time. If the graph models a network of interacting qubits, transfer of state among such qubits throughout time can be formalized as the action of the continuous-time quantum walk operator in the characteristic vectors of the vertices.

Here we are concerned with the problem of determining which graphs admit a perfect transfer of state. More specifically, we will study graphs whose adjacency matrix is a sum of tensor products of $01$-matrices, focusing on the case where a graph is the tensor product of two other graphs. As a result, we will construct many new examples of perfect state transfer.
\end{abstract}

\section{Introduction} 

We model a network of interacting qubits by a simple and undirected graph, and we address the problem of when such a network admits a transfer of quantum state between two sites without a loss of information. We call this phenomenon \textsl{perfect state transfer}.

Perfect state transfer is a relatively rare phenomenon. In \cite{GodsilPerfectStateTransfer12}, Godsil proved that there are only finitely many connected graphs with maximum valency $k$ in which perfect state transfer can occur. 

The main achievement of this paper is the construction of many new examples of simple graphs admitting perfect state transfer. These examples are obtained by examining certain graph products. In fact, we introduce a general framework in which perfect state transfer in graph products can be studied.

To the best of our knowledge, the following list contains a summary of which simple graphs are known to admit perfect state transfer. To highlight the importance of our findings, we point out that we are significantly increasing the quantity of known examples of perfect state transfer.

\begin{enumerate}
\item Paths $P_2$ and $P_3$ between vertices of degree $1$, but no other $P_n$ for $n\geq 4$ (see \cite{ChristandlPSTQuantumSpinNet}).
\item Certain cubelike graphs, fully characterized in \cite{GodsilBernasconiPSTCubelike} and \cite{GodsilCheungPSTCubelike}.
\item Certain circulant graphs, fully characterized in \cite{BasicPetkovicPSTCirculant} and \cite{BasicCirculant}.
\item Certain graphs admitting an equitable partition whose quotient is a weighted path. Among those, all hypercubes and cartesian powers of $P_3$ (see \cite{ChristandlPSTQuantumSpinNet}) and some distance-regular graphs (see preprint in \cite{CoutinhoGodsilGuoVanhove}).
\item Bipartite doubles of certain graphs belonging to association schemes (see preprint in \cite{CoutinhoGodsilGuoVanhove}).
\item Joins of the form $K_2 + G$ or $\overline{K_2} + G$, where $G$ is any regular graph with certain specified orders and degrees (see \cite{Angeles-CanulPSTcirculant}).
\item Iterated self-joins of regular graphs admitting perfect state transfer for certain choices of the parameters (see \cite{Angeles-CanulPSTcirculant}).
\item Certain circulant joins of circulant graphs (see \cite{Angeles-CanulPSTcirculant}).
\item Graphs obtained from lifting on certain weighted paths on four vertices (see \cite{BachmanTamonPSTQuotientandLifting}).
\item Some graphs which are the products of other graphs. Determined in \cite{ChristandlPSTQuantumSpinNet2} for the cartesian product. Some examples involving other products of graphs were constructed in \cite{GeGreenbergPerezTamonPSTproducts} and \cite{NEPSP3}.
\end{enumerate}

In this paper, we almost fully characterize when a tensor product of two graphs $X$ and $Y$ admits perfect state transfer. As a corollary of our work, will show that under certain mild integrality conditions on the eigenvalues of the graphs, if $X$ admits perfect state transfer and the power of two in the factorization of the integer parts of the eigenvalues of $Y$ is constant, then $X^{\square k} \times Y$ admits perfect state transfer for an infinite number of values of $k$. For instance, $Y$ can be a star or any graph with odd integer eigenvalues.

We will also apply our methods to determine when switching graphs admit perfect state transfer. Here again we will construct a new infinite family of graphs admitting perfect state transfer.

All of the results in this paper and an elementary introduction to the topic of perfect state transfer can be found in \cite[PhD Thesis]{CoutinhoPhD}.

\section{Definitions and basic results} \label{section:1}

Let $X$ be a (simple and undirected) graph and consider its adjacency matrix $A = A(X)$. For every non-negative real number $t$, we denote ${U_A(t) = \exp(\ii t A)}$, thus
\[U_A(t)  = \sum_{k \geq 0} \frac{(\ii t)^k}{k!} A^k.\]
We omit the subscript $A$ whenever the context is clear. Because $A$ is symmetric, $U(t)$ is a unitary operator.

We say that a graph $X$ admits \textit{perfect state transfer} at time $\tau$ from vertex $u$ to a distinct vertex $v$ if
\[U(\tau) \ee_u = \lambda \ee_v\]
where $\ee_u$ and $\ee_v$ are the characteristic vectors of the vertices $u$ and $v$ and $\lambda$ is a complex number of absolute value $1$, called the \textit{phase}. Note that $|\lambda| = 1$ because $U(t)$ is unitary for all $t$. Note also that perfect state transfer from $u$ to $v$ implies perfect state transfer from $v$ to $u$. For this reason, we might refer to $uv$-perfect state transfer in a graph $X$. We also say that $X$ is \textit{periodic} at $u$ at time $\tau \neq 0$ if
\[U(\tau) \ee_u = \lambda \ee_u\]
for some $\lambda \in \mathbb C$. Note that if $X$ admits $uv$-perfect state transfer at time $\tau$, then $X$ is periodic at both $u$ and $v$ at time $2 \tau$.

Suppose that the adjacency matrix $A$ of a graph has distinct eigenvalues $\theta_0 > \theta_1 > ... > \theta_d$. The matrix $A$ is real and symmetric, so it is orthogonally diagonalizable and therefore admits a spectral decomposition
\[A = \sum_{r=0}^d \theta_r E_r\]
into a linear combination of orthogonal projectors. Given a graph $X$ and a vertex $u \in V(X)$ with characteristic vector $\ee_u$, the \textit{eigenvalue support} of $u$, to be denoted by $\Phi_u$, is defined as the set of eigenvalues of $X$ such that the projection of $\ee_u$ into the corresponding eigenspace is a non-zero vector. Let $v$ be another vertex. If, for all $r \in \{0,...,d\}$, it follows that $E_r \ee_u = \pm E_r \ee_v$, we say that $u$ and $v$ are strongly cospectral vertices. Note that if $u$ and $v$ are strongly cospectral, then $\Phi_u = \Phi_v$, and in this case, let us define $\Phi_{uv}^+$ as the set of eigenvalues $\theta_r$ in $\Phi_u$ and such that $E_r \ee_u = E_r \ee_v$, and $\Phi_{uv}^-$ as the set of eigenvalues $\theta_r$ in $\Phi_u$ and such that $E_r \ee_u = - E_r \ee_v$.

The set up above allows us to translate the phenomenon of perfect state transfer into conditions that depend uniquely on the spectral decomposition of $A$. Note that even though the ideas in the theorem below have been used in previous works to characterize perfect state transfer in various classes graphs, we do not know of any published paper in which this theorem appears with this level of generality. For a detailed proof, please check \cite{CoutinhoPhD}.

\begin{theorem} \label{thm:2}
Let $X$ be a graph with vertices $u$ and $v$ and distinct eigenvalues $\theta_0 > ... > \theta_d$. There exists a time $\tau \in \mathbb{R^+}$ and a constant $\lambda \in \mathbb{C}$ satisfying
\[U_X(\tau) \ee_u = \lambda \ee_v\]
if and only if the following conditions hold.
\begin{enumerate}[(i)]
\item Vertices $u$ and $v$ are strongly cospectral.
\item \cite[Theorem 6.1]{GodsilPerfectStateTransfer12} Non-zero elements in $\Phi_u$ are either all integers or all quadratic integers. Moreover, there is a square-free integer $\Delta$, an integer $a$, and integers $b_0,...,b_k$ such that
\[\theta_r = \frac{1}{2} \big( a + b_r \sqrt{\Delta} \big) \ \ \text{for all $r = 0,...,k$}.\]
Here we allow $\Delta = 1$ for the case where all eigenvalues are integers, and $a =0$ for the case where they are all multiples of $\sqrt{\Delta}$.
\item \label{pst:thm:pstmain2(iii)} Let $\displaystyle g = \gcd\left( \left\{ \frac{\theta_0 - \theta_r}{\sqrt \Delta} \right\}_{r = 0}^k \right)$. Then 
\begin{enumerate}[a)]
\item $\theta_r \in \Phi_{uv}^+$ if and only if $\dfrac{\theta_0 - \theta_r}{g \sqrt{\Delta}}$ is even, and
\item $\theta_r \in \Phi_{uv}^-$ if and only if $\dfrac{\theta_0 - \theta_r}{g \sqrt{\Delta}}$ is odd.
\end{enumerate}
\end{enumerate}
Moreover, if these conditions hold and perfect state transfer occurs between $u$ and $v$ at time $\tau$ with phase $\lambda$, then
\begin{enumerate}[a)]
\item There is a minimum time $\tau_0>0$ at which perfect state transfer occurs between $u$ and $v$, and
\[\tau_0 = \frac{1}{g} \frac{\pi}{\sqrt{\Delta}}.\]
\item The time $\tau$ is an odd multiple of $\tau_0$.
\item The phase $\lambda$ is equal to $\e^{\ii \tau \theta_0}$. \qed
\end{enumerate} 
\end{theorem}

\section{Graph products} \label{section:3}
\newcommand{\properties}{\textsc{properties}}
\newcommand{\A}{\mathcal{A}}
\newcommand{\B}{\mathcal{B}}

In this section, we briefly introduce a very general method for solving the problem of finding perfect state transfer in graphs whose adjacency matrices are sums of tensor products of $01$-matrices. For the sake of simplifying the notation, we will restrict our considerations to graphs $X$ such that
\[A(X) = B \otimes C + M \otimes N\]
with $B,C,M,N$ symmetric matrices. This includes most graph products, but both the number of terms of the sum and the number of factors in each term can be generalized to any positive integer. Example of graphs of this form are the cartesian product $X \square Y$ of graphs $X$ and $Y$, whose adjacency matrix satisfies 
\[A(X \square Y) = A(X) \otimes I + I \otimes A(Y).\] 
The existence of perfect state transfer in the cartesian product is easily determined from the following equation (see \cite{ChristandlPSTQuantumSpinNet2}).
\begin{align}U_{X \square Y}(t) = U_X(t) \otimes U_Y(t). \label{eq:5} \end{align}

\begin{lemma} \label{lem:4}
Suppose $X$ is a graph such that $A(X) = B \otimes C + M \otimes N$, where $B$ and $M$ are commuting $m \times m$ matrices. Let $\beta_0 \geq ...\geq \beta_{m-1}$ and $\mu_0 \geq ... \geq \mu_{m-1}$ be the spectra of $B$ and $M$ respectively. Then $U_{X}(t)$ is similar to a block diagonal matrix with $m$ blocks, in which the $r$-th block is equal to $\exp(\ii t L_r)$, with \[L_r = \beta_r C + \mu_r N.\]
\end{lemma}

\proof
Because $B$ and $M$ are symmetric and commute, they can be simultaneously diagonalized by a matrix $P$. As a consequence, $A(X)$ is similar to a block diagonal matrix whose blocks are equal to $L_r$, with $r \in \{0,...,m-1\}$. The result now follows trivially.
\qed

A typical vertex of $X$ will be represented as $(w,u)$, where $w$ indexes a row of $B$ and $M$, and $u$ a row of $C$ and $N$. In the context of the lemma above, we allow the use of the terms \textsl{perfect state transfer} and \textsl{periodicity} to symmetric matrices in general.

\begin{lemma}\label{lem:5}
Suppose $X$ is defined as in Lemma \ref{lem:4}, and let $E_0,...,E_{m-1}$ be the rank-1 projectors onto the common eigenlines of $B$ and $M$. Let $\displaystyle L_r = \beta_r C + \mu_r N$. If the graph $X$ admits perfect state transfer from $(w,u)$ to $(z,v)$, then
\begin{enumerate}[(i)]
\item the vertices $w$ and $z$ are strongly cospectral in the matrices $B$ and $M$; and
\item for all $r$ such that $E_r w \neq 0$; if $u = v$, the matrix $L_r$ is periodic at $u$; and if $u \neq v$, the matrix $L_r$ admits perfect state transfer from $u$ to $v$.
\end{enumerate}
\end{lemma}
\proof
By hypothesis, there exists a time $\tau$ and a complex number $\lambda$ such that
\[U_{X}(\tau) (\ee_w \otimes \ee_u) = \lambda (\ee_z \otimes \ee_v).\]
Let $P$ be the matrix that simultaneously diagonalizes $B$ and $M$. Thus
\[(P^T \otimes I) U_{X}(\tau) (P \otimes I) ( P^T \ee_w \otimes \ee_u) = \lambda (P^T \ee_z \otimes \ee_v),\]
where $(P^T \otimes I) U_{X}(\tau) (P \otimes I)$ is the block diagonal matrix of Lemma \ref{lem:4}. This is true if and only if, for all $r = 0,...,m-1$,
\[(\ee_r^T P^T \ee_w) \exp(\ii \tau L_r) \ee_u = (\ee_r^T P^T \ee_z)  \lambda \ee_v.\]
Because the matrices $\exp(\ii \tau L_r)$ are all unitary, it follows that, for all $r$,
\[(\ee_r^T P^T \ee_w) = \pm (\ee_r^T P^T \ee_z),\]
and so $w$ and $z$ are strongly cospectral in the matrices $B$ and $M$; and also that
\[\exp(\ii \tau L_r)  \ee_u = \pm \lambda \ee_v\]
for all $r$ such that $(\ee_r^T P^T \ee_w)  \neq 0$.
\qed

\begin{lemma}\label{lem:generalizedprods}
Suppose $X$ is defined as in Lemma \ref{lem:5}, and let $\gamma_0 \geq ... \geq \gamma_{n-1}$ and $\nu_0 \geq ... \geq \nu_{n-1}$ be the spectra of $C$ and $N$ respectively. Suppose that $C$ and $N$ also commute, and let $F_0,...,F_{n-1}$ be the rank-1 projectors onto the common eigenlines of $C$ and $N$. Then
\[U_{X}(t) =\sum_{r = 0}^{m-1} \sum_{s=0}^{n-1} \e^{\ii t(\beta_r\gamma_s + \mu_r\nu_s)} E_r \otimes F_s .\]
\end{lemma}

\proof
Note that
\[B\otimes C+ M \otimes N= \sum_{r=0}^{m-1}\sum_{s=0}^{n-1} (\beta_r\gamma_s + \mu_r\nu_s) E_r \otimes F_s.\]
From
\[U_{X}(t) = \exp\Big(\ii t \big(B\otimes C+ M\otimes N\big)\Big),\]
and the facts that $BM = MB$ and $CN = NC$, it follows that
\[U_{X}(t) =\sum_{r = 0}^{m-1} \sum_{s=0}^{n-1} \e^{\ii t (\beta_r\gamma_s + \mu_r\nu_s)  } E_r \otimes F_s.\]
\qed

\section{Tensor product of graphs} \label{section 2}

Consider graphs $X$ and $Y$ with respective adjacency matrices $A(X)$ and $A(Y)$. The \textit{tensor product} of $X$ and $Y$, denoted by $X \times Y$, is the graph with adjacency matrix $A(X) \otimes A(Y)$. This product is also known in the literature as the weak direct product, the direct product, the categorical product, and many other names (see \cite[Chapter 4]{GraphProds}). 

In this section, we almost completely determine when a tensor product of graphs admits perfect state transfer. As an immediate application, we find more examples of perfect state transfer.

\begin{theorem} \label{thm:1}
Suppose $X$ and $Y$ are graphs, and $X \times Y$ admits perfect state transfer between vertices $(w,u)$ and $(z,v)$. If $u = v$, then $Y$ is periodic at $u$. If $u \neq v$, then $Y$ admits $uv$-perfect state transfer. Likewise, if $w = z$, then $X$ is periodic at $w$. If $w \neq z$, then $X$ admits $wz$-perfect state transfer.
\end{theorem}
\proof
It is a simple consequence of Lemma \ref{lem:5} with $B=C=0$, $M=A(X)$ and $N = A(Y)$, or $M=A(Y)$ and $N = A(X)$.
\qed

\begin{lemma} \label{lem:1}
Suppose $X$ and $Y$ are graphs and $A(X)$ admits the spectral decomposition $\displaystyle A(X) = \sum_{r = 0}^d \theta_r E_r$. Then
\[U_{X \times Y} (t) = \sum_{r=0}^{d} E_r \otimes U_Y(\theta_r t).\]
\end{lemma}
\proof
It is a consequence of Lemma \ref{lem:generalizedprods} and an easy rearrangement.
\qed

Now we will determine under which conditions on the factors we can obtain perfect state transfer on the product.
\newcommand{\s}{m}
\begin{theorem} \label{thm:3}
Suppose $U_Y(\tau)\ee_u = \lambda \ee_v$, and that the eigenvalues of $Y$ in the support of $u$ are of the form $b_i \sqrt{\Delta_u}$. Suppose $w$ and $z$ are strongly cospectral vertices in $X$. Then $X \times Y$ admits perfect state transfer from $(w ,u)$ to $(z ,v)$ if and only if the following conditions hold.
\begin{enumerate}[(i)]
\item For all $\theta_r \in \Phi_w$, we have $\theta_r = t_r \sqrt{\Delta_w}$, where $t_r \in \mathbb{Z}$ and $\Delta_w$ is a square-free positive integer (which could be $1$).
\item The powers of two in the factorizations of each $t_r$ are all equal.
\item If $\lambda$ is a primitive $n$-th root of the unit, then $n$ is even, and there exists an integer $\s$ such that
\begin{enumerate}[a)]
\item If $\theta_r \in \Phi_{wz}^+$, then the odd part of $t_r$ is congruent to $\s$ modulo $n$.
\item If $\theta_r \in \Phi_{wz}^-$, then the odd part of $t_r$ is congruent to $\s + \frac{n}{2}$ modulo $n$.
\end{enumerate}
\end{enumerate}
\end{theorem}

\proof
Let $\Phi_w = \{\theta_0,...,\theta_d\}$, and $\Phi_u = \{\varphi_0,...,\varphi_k\}$. Let $h$ be the gcd of the differences $(\varphi_0-\varphi_r)$ for $r=1,...,k$. Let $h = 2^e \ell$, with $\ell$ an odd integer.

Suppose that perfect state transfer occurs in $X \times Y$ between $(w,u)$ and $(z,v)$ at time $\tau$ and phase $\gamma$. As a consequence of the fact that $A(X \times Y) = A(X) \otimes A(Y)$, the eigenvalues in the support of $(w,u)$ are of the form $\theta_r \varphi_i $, with $0 \leq r \leq d$ and $0 \leq i\leq  k$. In light of Theorem \ref{thm:2}, the eigenvalues $\theta_r$ are either integers or integer multiples of $\sqrt{\Delta_w}$ for some square-free positive $\Delta_w \in \mathbb{Z}$.

Then, using Lemma \ref{lem:1}, we have
\begin{align*}
\gamma (\ee_z \otimes \ee_v) & = U_{X \otimes Y} (\tau) (\ee_w \otimes \ee_u) \\ & = \sum_{r=0}^d (E_r \otimes U_Y(\theta_r \tau) )(\ee_w \otimes \ee_u) \\ 
& = \sum_{r=0}^d E_r \ee_w \otimes U_Y(\theta_r \tau) \ee_u.
\end{align*}
Multiplying both sides by $E_r \otimes I$, we get that, for $\sigma_r \in \{+1,-1\}$,
\[U_Y(\theta_r \tau) \ee_u = \sigma_r \gamma \ee_v\]
depending on whether $\theta_r \in \Phi_{wz}^+$ or $\theta_r \in \Phi_{wz}^-$. In either case, $\theta_r \tau$ is a time for which perfect state transfer occurs in $Y$ between $u$ and $v$. Applying Theorem \ref{thm:2}, this implies that
\[\theta_r \tau = \ell_r \frac{\pi }{2^e . \ell \sqrt{\Delta_u}} \quad \text{and} \quad \sigma_r \gamma = \lambda^{\ell_r},\]
where $\ell_r$ is an odd integer. Considering $\theta_r$ and $\theta_s$ in the support of $w$, we will have
\begin{align} \label{eq:4} \frac{\theta_r}{\theta_s} = \frac{\ell_r}{\ell_s}.\end{align}
Because the integers $\ell_r$ are odd, the power of $2$ dividing each $t_r$ is the same, proving condition (ii).

To prove condition (iii), suppose we take $\s' \in \{1,...,n\}$ such that $\lambda^{\s'} = \gamma$. The fact that there is a $\s''$ such that $\lambda^{\s''} = - \gamma$ is equivalent to $(-1)$ being a power of $\lambda$, which happens if and only if $n$ is even. In that case, if $\theta_r \in \Phi_{wz}^+$, then $\ell_r \equiv \s' \bmod n$, and if $\theta_r \in \Phi_{wz}^-$, then $\ell_r \equiv \s' + \frac{n}{2} \bmod n$. Note that the odd part of $t_r$ is an odd multiple $\ell_r$, which by Equation \ref{eq:4} does not depend on $r$. Say $\ell'$. So if the integers $\ell_r$ satisfy the congruences with $\s'$, so will the odd parts of $t_r$ with $\s = \s' \ell'$.

Now suppose all three conditions hold. Let $t_r = 2^f.k_r$, for some $f \geq 0$ and odd integers $k_r$. By Lemma \ref{lem:1}, we have
\[U_{X \times Y} \left( \frac{\pi }{2^{e+f} \ell \sqrt{\Delta_w}\sqrt{\Delta_u}} \right) = \sum_{r=0}^d E_r \otimes U_Y\left( \theta_r . \frac{\pi }{2^{e+f}  \ell \sqrt{\Delta_w}\sqrt{\Delta_u}} \right).\]
Note that
\begin{align*}
U_Y\left( \theta_r . \frac{\pi }{2^{e+f}  \ell \sqrt{\Delta_w}\sqrt{\Delta_u}} \right) &  = \left[U_Y\left( \frac{\pi}{2^e \ell \sqrt{\Delta_u}}\right)\right]^{k_r }
\end{align*}
Hence
\begin{align*}
U_{X \times Y} \left( \frac{\pi }{2^{e+f}  \ell \sqrt{\Delta_w}\sqrt{\Delta_u}} \right) (\ee_w \otimes \ee_u) & = \sum_{r=0}^d E_r \ee_w \otimes \left[U_Y\left( \frac{\pi}{2^e\sqrt{\Delta_u}}\right)\right]^{k_r } \ee_u
\\ & = \sum_{r=0}^d E_r \ee_w \otimes \lambda^{k_r} \ee_v.
\end{align*}
If $r \in \Phi_{wz}^+$, then condition (iii) implies that $\lambda^{k_r} = \lambda^{\s}$, and if $r \in \Phi_{wz}^-$, then $\lambda^{k_r} = \lambda^{-\s}$. If $\lambda^{\s} = \gamma$, we have
\[U_{X \times Y} \left( \frac{\pi }{2^{e+f}  \sqrt{\Delta_w}\sqrt{\Delta_u}} \right) (\ee_w \otimes \ee_u) = \gamma (\ee_z \otimes \ee_v),\]
as claimed.
\qed

Note that if $w=z$, the theorem above takes the following form.

\begin{corollary} \label{cor:1}
Suppose $U_Y(\tau)\ee_u = \lambda \ee_v$, and that the eigenvalues of $Y$ in the support of $u$ are of the form $b_i \sqrt{\Delta_u}$. Let $w \in V(X)$. Then $X \times Y$ admits perfect state transfer from $(w ,u)$ to $(w ,v)$ if and only if the following conditions hold.
\begin{enumerate}[(i)]
\item For all $\theta_r \in \Phi_w$, we have $\theta_r = t_r \sqrt{\Delta_w}$, where $t_r \in \mathbb{Z}$ and $\Delta_w$ is a square-free positive integer (which could be $1$).
\item The powers of two in the factorizations of each $t_r$ are all equal.
\item If $\lambda$ is a primitive $n$-th root of the unit, then there exists an integer $\s$ such that the odd part of the integer $t_r$ is congruent to $\s$ modulo $n$.
\end{enumerate}
\end{corollary}

%TODO remainig case for eigenvalues: a+b_i\sqrtdelta

Note that the conditions on $X$ of both results depend very little on $Y$. In fact, if $\varphi_0$ is the largest eigenvalue of $Y$ and $\tau$ is the time at which perfect state transfer occurs in $Y$, then the conditions depend only on the eigenvalues of $X$, except for the order of $\e^{\ii  \varphi_0 \tau}$ as a root of unity. As a consequence, we have the following result. We will use $Y^{\square k}$ to denote the Cartesian product of the graph $Y$ with itself $k$ times.

\begin{corollary}
If $X \times Y$ admits perfect state transfer,  and if the eigenvalues of $Y$ in the support of vertices involved in perfect state transfer are integers or integer multiples of a square root, then $X \times Y^{\square k}$ admits perfect state transfer for all $k \in \mathbb{Z}^+$.
\end{corollary}
\proof
By Equation (\ref{eq:5}), if $Y$ admits perfect state transfer at minimum time $\tau$, then so does $Y^{\square k}$. The corollary now follows from the fact that the largest eigenvalue of $Y^{\square k}$ is $k\varphi_0$, so the order of the phase as a root of unity does not increase.
\qed

This can be pushed even further.

\begin{corollary} \label{cor:nice1}
If $Y$ admits perfect state transfer, if the eigenvalues of $X$ and $Y$ are integers or integer multiples of a square root, and if the powers of two in the factorizations of the integer parts of the eigenvalues of $X$ are all the same, then there exists a $k_0 \in \mathbb{Z}^+$ such that $X \times Y^{\square (mk_0)}$ admits perfect state transfer for all $m \geq 1$.
\end{corollary}

As a consequence, we present new examples of perfect state transfer in simple graphs. For the cases below, we assume that $Y$ is a graph admitting $uv$-perfect state transfer, and the eigenvalues of $Y$ in the support of $u$ are either integers or integer multiples of square roots. All the graphs known in the literature admitting perfect state transfer are of this form.

\paragraph{Stars}%TODO CHECK WHEELS AND WINDMILLS

Let $S_n$ represent the graph on $n+1$ vertices with degree sequence $(n,1,1,1,...,1)$. The spectrum of $S_n$ is
\[ \{ \sqrt{n}^{(1)},\ 0^{(n-2)},\ -\sqrt{n}^{(1)}\}.\]
Let $w$ be the vertex of degree $n$. The eigenvalue support of $w$ is $\{\sqrt{n},-\sqrt{n}\}$. 

From Corollary \ref{cor:nice1}, there is a $k$ such that $S_n \times Y^{\square k}$ admits perfect state transfer from $(w,u,u,...,u)$ to $(w,v,v,...,v)$. 

Note that $k$ is usually quite small. If $Y$ admits perfect state transfer at time $\frac{\pi}{2}$, which is a rather common situation, then $k = 2$ will suffice.

\paragraph{Odd eigenvalues}

If $X$ is a graph with odd eigenvalues, and $w \in V(X)$, then it follows from Corollary \ref{cor:nice1} that there is a $k$ such that $X \times Y^{\square k}$ admits perfect state transfer from $(w,u,u,...,u)$ to $(w,v,v,...,v)$.

We can find many graphs with odd eigenvalues among the known distance-regular graphs. For example, there are 32548 non-isomorphic strongly regular graphs with parameters (36,15,6,6). These graphs have distinct eigenvalues $\{15,3,-3\}$. The tensor product of each of them with $C_4$ will admit perfect state transfer.

\section{Double covers and switching graphs} \label{section:switching}

Given graphs $X$ and $Y$ on the same set of vertices, we define the graph $X \ltimes Y$ as the graph with adjacency matrix
\[A(X \ltimes Y) = \begin{pmatrix} A(X) & A(Y) \\ A(Y) & A(X) \end{pmatrix}\]
If $A(X)$ and $A(Y)$ have their ones in disjoint positions, then $X \ltimes Y$ is a \textit{double cover} of the graph with adjacency matrix $A(X)+A(Y)$. When $A(Y) = J - I - A(X)$, $X \ltimes Y$ is a double cover of the complete graph and is known in the literature as the \textit{switching graph} of $X$ (see \cite[chapter 11]{GodsilRoyle}).

If $X$ is the empty graph, then $X \ltimes Y = A(K_2) \otimes A(Y)$ is known as the \textit{bipartite double} of $Y$. Perfect state transfer on bipartite doubles of graphs belonging to association schemes was studied in \cite{CoutinhoGodsilGuoVanhove}. Here we intend to study $X \ltimes Y$ in a more general form.

Note that $I_2$ and $A(K_2)$ commute, and can be simultaneously diagonalized by $H = \frac{1}{\sqrt{2}} \left(\begin{array}{cc} 1 & 1 \\ 1 & -1\end{array} \right) $. Using Lemma \ref{lem:4}, we have:

\begin{lemma} \label{lem:switching2}
Given graphs $X$ and $Y$, $A = A(X)$ and $B=A(Y)$, we have, for all $t \geq 0$,
\[(H \otimes I) U_{X \ltimes Y}(t)(H \otimes I)= \left( \begin{array}{cc}
U_{A+B}(t) & 0 \\ 0 & U_{A-B}(t)
\end{array}\right).\]
\end{lemma}

\begin{theorem} \label{thm:22}
Given graphs $X$ and $Y$ on the same vertex set $V$, with $A = A(X)$ and $B=A(Y)$, the graph $X \ltimes Y$ on vertex set $\{0,1\} \times V$ admits perfect state transfer if and only if one of the following holds.
\begin{enumerate}[(i)]
\item For some $\tau \in \mathbb{R}^+$, $\lambda \in \mathbb{C}$ and $u \in V$, the matrices $A+B$ and $A-B$ are periodic at $u$ at time $\tau$ with respective phase factors $+\lambda$ and $-\lambda$. In this case, perfect state transfer is between $(0,u)$ and $(1,u)$.
\item For some $\tau \in \mathbb{R}^+$, $\lambda \in \mathbb{C}$ and $u ,v\in V$, the matrices $A+B$ and $A-B$ admit $uv$-perfect state transfer at time $\tau$ with the same phase factor $\lambda$. In this case, perfect state transfer is between $(0,u)$ and $(0,v)$, and between $(1,u)$ and $(1,v)$.
\item For some $\tau \in \mathbb{R}^+$, $\lambda \in \mathbb{C}$ and $u ,v\in V$, the matrices $A+B$ and $A-B$ admit $uv$-perfect state transfer at time $\tau$ with respective phase factors $+\lambda$ and $-\lambda$. In this case, perfect state transfer is between $(0,u)$ and $(1,v)$, and between $(1,u)$ and $(0,v)$.
\end{enumerate}
\end{theorem}
\proof
From Lemma \ref{lem:switching2}, it follows that
\[U_{X \ltimes Y}(t)= \frac{1}{2} \left( \begin{array}{cc}
U_{A+B}(t) +  U_{A-B}(t) & U_{A+B}(t) -  U_{A-B}(t) \\ U_{A+B}(t) - U_{A-B}(t) &U_{A+B}(t) +  U_{A-B}(t)
\end{array}\right).\]
For any $u \in V$, perfect state transfer happens at $\tau$ between $(0,u)$ and $(1,u)$ if and only if
\[\big| \big(U_{A+B}(\tau) - U_{A-B}(\tau) \big)_{u,u} \big| = 2.\]
Because $|U_{A+B}(\tau)_{u,u}| \leq 1$ and $|U_{A-B}(\tau)_{u,u}| \leq 1$, this is equivalent to
\[\lambda = U_{A+B}(\tau)_{u,u} = -\ U_{A-B}(\tau)_{u,u} \quad \text{and} \quad |\lambda| = 1,\]
proving (i).

Likewise, perfect state transfer happens at $\tau$ between $(0,u)$ and $(0,v)$ if and only if
\[\big| \big(U_{A+B}(\tau) + U_{A-B}(\tau) \big)_{u,v} \big| = 2.\]
Because $|U_{A+B}(\tau)_{u,v}| \leq 1$ and $|U_{A-B}(\tau)_{u,v}| \leq 1$, this is equivalent to
\[\lambda = U_{A+B}(\tau)_{u,v} = U_{A-B}(\tau)_{u,v} \quad \text{and} \quad |\lambda| = 1.\]
Naturally these conditions are also equivalent to perfect state transfer at time $\tau$ between $(1,u)$ and $(1,v)$. This proves (ii).

Case (iii) can be proved similarly, and the three cases together cover all possible ways in which a vertex of $X \ltimes Y$ could be involved in perfect state transfer.
\qed

We show below how to obtain new examples of perfect state transfer in switching graphs using case (i).

\begin{corollary}
Let $\theta_0 > ... > \theta_d$ be the eigenvalues of $A(X) - A(\overline{X})$. Suppose $|V(X)|=n >2$. Then $X \ltimes \overline{X}$ admits perfect state transfer if and only if there is an integer $k$ such that, for all $\theta_r$ in the support of a vertex $u$, the number $2k(\theta_r + 1)/n$ is odd for all $r$.
\end{corollary}
\proof
Let $A = A(X)$ and $B = A(\overline{X})$. First note that $A+B = J-I$ is the adjacency matrix of the complete graph, thus the matrix $A+B$ does not admit perfect state transfer when $n>2$. So if $X \ltimes \overline{X}$ admits perfect state transfer, we must be in case (i) of Theorem \ref{thm:22}.

Let $\lambda = [U_{A+B}(\tau)]_{u,u}$. Because $A+B = J-I$, it follows that $|\lambda| = 1$ if and only if $\tau = \frac{2 k \pi}{n}$ for any integer $k$. In that case, $\lambda = \e^{-\ii  \tau}$. Given Theorem \ref{thm:22}, perfect state transfer between $(0,u)$ and $(1,u)$ happens if and only if
\[\e^{\ii  \tau \theta_r} = - \e^{-\ii \tau}\]
for all $\theta_r$ in the support of $u$. This is equivalent to the condition in the statement.
\qed

When $X$ is a regular graph, then $X$ and $\overline{X}$ commute. More generally, if $A=A(X)$ and $B = A(Y)$ commute, then $U_{A\pm B}(t) = U_A(t)U_B(\pm t)$. Using this fact, it is very easy to find examples of graphs $X$ and $Y$ such that $X \ltimes Y$ admits perfect state transfer using Theorem \ref{thm:22}.

For instance, if $X$ is the empty graph and $Y$ is a non-bipartite connected graph admitting perfect state transfer with phase $\pm \ii$, then $X \ltimes Y$ is connected and admits perfect state transfer. Similarly, if $X$ is a graph admitting perfect state transfer and $Y$ is such that $A(Y) = I$, then $X \ltimes Y$ is a simple graph admitting perfect state transfer. 

More interestingly, if $A(X)$ and $A(Y)$ commute, and if there is a time $\tau$ such that at this time $X$ is periodic at $u$ and $Y$ is periodic at $u$ with phase $\pm \ii$, then perfect state transfer occurs in $X \ltimes Y$ between $(0,u)$ and $(1,u)$. We use this to generate the examples below.

\paragraph{Switching graphs}

The graph $K_n \square K_n$ is a strongly regular graph with parameters $(n^2,2n-2,n-2,2)$. Its eigenvalues are $\{2n-2, \ n-2,\ -2 \}$, hence for all $n$ divisible by $4$, it follows from the observation above that the switching graph of $K_n \square K_n$ admits perfect state transfer at time $\pi / 2$.

There are two feasible parameter sets for strongly regular graphs on 96 vertices for which constructions of such graphs are known (see \cite{BrouwerKoolenKlin}). The parameter sets are $(96,20,4,4)$ and $(96,76,60,60)$. One example for the first set is the point-graph of the generalized quadrangle GQ(5,3). The switching graph of all strongly regular graphs with such parameters admit perfect state transfer.

\section*{Acknowledgements}

We thank an anonymous referee for some invaluable remarks.

\bibliographystyle{plain}
\bibliography{qwalks}

\end{document}